\documentclass{article}

\usepackage{arxiv}

\usepackage[utf8]{inputenc} % allow utf-8 input
\usepackage[T1]{fontenc}    % use 8-bit T1 fonts
\usepackage[hidelinks]{hyperref}       % hyperlinks
\usepackage{url}            % simple URL typesetting
\usepackage{booktabs}       % professional-quality tables
\usepackage{amsfonts}       % blackboard math symbols
\usepackage{nicefrac}       % compact symbols for 1/2, etc.
\usepackage{microtype}      % microtypography
\usepackage{lipsum}         % Can be removed after putting your text content
\usepackage{graphicx}
\usepackage{doi}

\usepackage{colortbl}
\usepackage{listings}
\usepackage{bm}

\usepackage{amsmath}
\usepackage{cleveref}       % smart cross-referencing

\usepackage[ruled,vlined]{algorithm2e}

\usepackage{appendix}
\usepackage{cite}

\usepackage{multirow}
\usepackage{array}
\newcommand{\head}[2]{\multicolumn{1}{>{\arraybackslash}p{#1}}{{#2}}}

\usepackage{chngcntr}
\SetKwComment{Comment}{/* }{ */}

\title{Separable Physics-Informed Neural Networks for the Solution of Elasticity Problems}

\date{\today}

\author{ \href{https://orcid.org/0000-0002-4930-1846}{\includegraphics[scale=0.06]{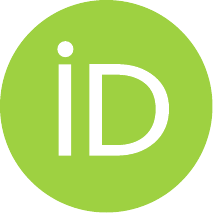}\hspace{1mm}Vasiliy A. Es'kin}\thanks{Corresponding author: Vasiliy Alekseevich Es’kin (vasiliy.eskin@gmail.com)} \\
	Department of Radiophysics, University of Nizhny Novgorod\\
	 Nizhny Novgorod, Russia, 603950\\
	 and\\
	 Huawei Nizhny Novgorod Research Center\\
	 Nizhny Novgorod, Russia\\
	\texttt{vasiliy.eskin@gmail.com} \\
	\And
	{\hspace{1mm}Danil V. Davydov} \\
	Mechanical Engineering Research Institute \\ 
	Russian Academy of Sciences\\
	Nizhny Novgorod, Russia, 603155\\
	and\\
	Huawei Nizhny Novgorod Research Center\\
	Nizhny Novgorod, Russia\\
	\texttt{davidovdan27@yandex.ru} \\
	\And
	{\hspace{1mm} Julia V. Gur'eva} \\
	Huawei Nizhny Novgorod Research Center\\
	Nizhny Novgorod, Russia\\
	{\texttt{gureva-yulya@list.ru}}
	\texttt{} \\
	\And
	\href{https://orcid.org/0000-0001-5210-8281}{\includegraphics[scale=0.06]{orcid.pdf}\hspace{1mm}Alexey O. Malkhanov} \\
	Huawei Nizhny Novgorod Research Center\\
	Nizhny Novgorod, Russia\\	\texttt{alexey.malkhanov@gmail.com} \\
	\And
	\href{https://orcid.org/0000-0002-0454-5249}{\includegraphics[scale=0.06]{orcid.pdf}\hspace{1mm}Mikhail E. Smorkalov} \\
	Skolkovo Institute of Science and Technology\\
	Moscow, Russia\\
	and\\
	Huawei Nizhny Novgorod Research Center\\
	Nizhny Novgorod, Russia\\	\texttt{smorkalovme@gmail.com} \\
}

\renewcommand{\headeright}{}
\renewcommand{\undertitle}{}
\renewcommand{\shorttitle}{A preprint}

\renewcommand{\vec}{\bf}

\hypersetup{
pdftitle={A template for the arxiv style},
pdfsubject={math.na, cs.AI, physics.comp-ph},
pdfauthor={Vasiliy A. Es'kin}
pdfkeywords={Deep Learning, Physics-informed Neural Networks, Separable Physics-Informed Neural Networks, Partial differential equations, Deep Energy Method, Predictive modeling, Computational physics},
}

\begin{document}
\maketitle

\begin{abstract}
	A method for solving elasticity problems based on separable physics-informed neural networks (SPINN) in conjunction with the deep energy method (DEM) is presented. 
	Numerical experiments have been carried out for a number of problems showing that this method has a significantly higher convergence rate and accuracy than the vanilla physics-informed neural networks (PINN) and even SPINN based on a system of partial differential equations (PDEs). In addition, using the SPINN in the framework of DEM approach it is possible to solve problems of the linear theory of elasticity on complex geometries, which is unachievable with the help of PINNs in frames of partial differential equations. Considered problems are very close to the industrial problems in terms of geometry, loading, and material parameters.
\end{abstract}

% keywords can be removed

\section{Introduction}
The recent few years have been marked by revolutionary achievements in the field of machine learning referred to Deep Neural Networks (DNN) and Deep Learning (DL)~\cite{Alzubaidi2021}. Continued improvement of DNN architecture and refinement of the training methods have led to the various applications which possess the capabilities previously available only for humans~\cite{Lecun2015,Linardatos2021,Khan2022}. {Neural networks based approaches are widely spread in} image classification, handwriting recognition, speech recognition, and translation, text generators, game systems and computer vision~\cite{Vaswani2018,Wang2019,Yao2020,Ramesh2021,FROLOV2021187,Silver2018,Schrittwieser2020,Ouyang2022}. This revolutionary {transformation} and the implementation of {deep learning} achievements in many areas of professional activity has not left aside the scientific and engineering fields. Developed DNN systems for solving scientific and engineering problems can be separated into those that leverage previously established physical laws (PINN)~\cite{Raissi2019,Wang2022,Wang2022_2,eskin2023optimal,wang2023experts}, those that utilized previously obtained data (from computations or experiments) and build their predictions on this basis~\cite{Lu_2021,Li2022,Krinitskiy2022,Fanaskov2022,Ovadia2023,jin2023mechanical}, and those that {exploit} both these approaches~\cite{Raissi2017, Raissi2017_2,Cai_2021}. In our paper we exclusively focus on PINN methods, avoiding data-driven approaches.

The main ideas for the training of neural networks, which predict the physical system behaviour, using physical laws have been around for decades\cite{Psichogios1992AHN,Lagaris_1998}. 
Due to modern computational tools, these ideas incarnated a few years ago in new approaches, such as (PINN)~\cite{Raissi2019}. The neural networks of these methods are trained to approximate the dependences of physical values {with the help of} equations together with initial and boundary conditions. This approach has been used to solve a wide range of problems described by ordinary differential equations~\cite{Raissi2017,Raissi2017_2,Raissi2019, Rackauckas2020}, integro-differential equations~\cite{Yuan2022}, nonlinear partial differential equations~\cite{Raissi2017,Raissi2017_2,Raissi2019,Jin2021,Zhao2021,Kharazmi2021}, PDE with noisy data~\cite{Yang2021}, etc.~\cite{Cuomo2022,Pang2020}, related to various fields of knowledge such as thermodynamics~\cite{Patel2022,LIU2024119565}, hydrodynamics~\cite{Cai2021,Jin2021,Zhao2021}, electrodynamics~\cite{Lin2021,Kharazmi2021}, geophysics~\cite{He2020}, finance~\cite{Yuan2022}, mechanics~\cite{HAGHIGHAT2021113741,MENG2023116172,Samaniego2020,NING2023115909,HAO2023107536,Moseley2020,NING2023116430} and thermo-mechanically coupled systems~\cite{Harandi_2023}.

Despite the PINN success in solid mechanics~\cite{HAGHIGHAT2021113741,MENG2023116172,Samaniego2020,NING2023115909,HAO2023107536}, for the 3D elasticity problem under complex geometry the PINN based only partial differential equations, which described displacements and stresses demonstrates low convergence. To address this drawback differential equations can be used along with conservation laws~\cite{Wang2022_2} in the PINN approach, or only the {energy conservation law}~\cite{Samaniego2020,NING2023116430}, which result in accuracy improvement for the considered problems. Another concern of the PINN is the low speed of training (but high speed of the evaluation) which can be solved by novel types of neural networks such as separable physics-informed neural networks (SPINN)~\cite{cho2023}. {The training time of SPINN is comparable with the time of obtaining an acceptable solution with the help of classical numerical methods like the finite element method (FEM).}

{This work is devoted to the application of a number of techniques that can improve the accuracy of prediction and the speed of training and evaluation on the basis of the PINN approach for the 3D elasticity problems. We will briefly list the contributions made in the paper}. We will briefly list the contributions made in the paper:
\begin{enumerate}
	\item The method for solving elasticity problems based on separable physics-informed neural networks in conjunction with the deep energy method is proposed.
	\item A comparison of the speed and accuracy of training for different PINNs is carried out: PINN based on PDE, SPINN based on PDE, SPINN based on a minimum of energy.
	\item It demonstrates that SPINN in the framework DEM approach solves problems of the linear elasticity on complex geometries, which is {unachievable with PINNs in terms} of partial differential equations.
\end{enumerate}

The paper is structured as follows. In section 2, the statement of the problem of linear elasticity on the basis of both the partial differential equation and energy minimization is presented. Section 3 includes the description of PINN approach for two different statements of the problem. Descriptions of architectures of the artificial neural networks, that were used in the work, are also provided here. In Section 4, numerical experiments to illustrate the accuracy and efficiency of the presented approach are given. Finally, in Section 5 concluding remarks are given.

\section{Statement of the Problem}

\subsection{Formulation of the problem on the basis of partial differential equations}

Consider a deformed body in a state of static equilibrium under the effect of external volumetric force ${\vec f}$ and external force on unit area of surface ${\vec T}$. Such an object is described by the following system of partial differential equations (PDE)~\cite{Landau1986_7}
\begin{align}\label{eq1}
	\frac{\partial \sigma_{ik}}{\partial x_k} + f_i = 0,\,\quad {\vec x}\in \Omega.  
\end{align}
and boundary conditions
\begin{align}
	\sigma_{ik}n_{k} = T_{i}, \quad {\vec x} \in \partial \Omega,\label{eq2_1}\\
	{\vec u} = {\vec u}^{(b)}, \quad {\vec x} \in \partial \Omega.\label{eq2_2}
\end{align}
Here, ${\vec x}$ are spatial coordinates, $\sigma_{ik}$ is the stress tensor,  ${\vec n}$ is a unit vector along the outward normal to the surface, ${\vec u}$ is displacement vector, ${\vec u}^{(b)}$ is displacement vector at boundaries, $\Omega$ and $\partial \Omega$ are the spatial domain of body and its boundary, respectively. The conditions (\ref{eq2_1}) and (\ref{eq2_2}) must be satisfied at every point of the surface of a body in equilibrium. Here and further, Einstein summation notation is used which dictates repeated indices should be summed. 

The stress tensor is given in terms of the strain tensor by {Hook law}
\begin{equation}\label{eq3}
	\sigma_{ik} = \lambda \varepsilon_{ll} \delta_{ik} +  \mu \varepsilon_{ik},
\end{equation}
where the strain tensor $\varepsilon_{ik}$ is
\begin{align}
	& \varepsilon_{ik} = \frac{1}{2} \left( \frac{\partial u_i}{\partial x_{k}} + \frac{\partial u_{k}}{\partial x_{i}}\right),\label{eq4}
\end{align}
the coefficients {(Lame constants)} are
\begin{align}
	& \lambda = \frac{\nu E}{2 (1+\nu)(1 - 2 \nu)},\quad \mu = \frac{E}{2(1+\nu)} \label{eq5},
\end{align}
$ \delta_{ik}$ is Kronecker symbol, $E$ and $\nu$ are Young's modulus and Poisson's ratio, respectively.

In all examples, we assume volumetric force is the gravitation force ${\vec f} = \rho {\vec g}$; $\rho$ is the density of body material and ${\vec g}$ is the gravitational acceleration vector.\\

In the Cartesian coordinate system ($x,y,z$) the system of equations (\ref{eq1}) has the following form
\begin{align}
	& \frac{\partial \sigma_{xx}}{\partial x} + \frac{\partial \sigma_{xy}}{\partial y} + \frac{\partial \sigma_{xz}}{\partial z} = 0,     \label{eq6}\\
	& \frac{\partial \sigma_{xy}}{\partial x} + \frac{\partial \sigma_{yy}}{\partial y} + \frac{\partial \sigma_{yz}}{\partial z} = 0, \label{eq7} \\
	& \frac{\partial \sigma_{xz}}{\partial x} + \frac{\partial \sigma_{yz}}{\partial y} + \frac{\partial \sigma_{zz}}{\partial z} - \rho g= 0.\label{eq8}	
\end{align}
The axes of a given Cartesian coordinate system are oriented so that the vector ${\vec g}$ is directed in the negative direction of $z$--axis. Here, the symmetrical properties of the stress tensor ($\sigma_{ik} = \sigma_{ki}$) are used.

Equations (\ref{eq1}) together with boundary conditions (\ref{eq2_1}) and (\ref{eq2_2}) are the total description of the steady-state linear elasticity problem which has a unique solution. This problem can be solved by using the classical finite-element method (FEM) or by a novel physics-informed neural network approach (PINN)~\cite{Raissi2017,Raissi2017_2,Raissi2019}.

\subsection{Formulation of the problem on the basis of energy minimization}

Instead of the PDE (\ref{eq1}) with boundary conditions (\ref{eq2_1}) and (\ref{eq2_2}) the steady-state linear elasticity problem can be considered in the following energy approach \cite{Landau1986_7,Samaniego2020}	
\begin{equation}
	{\vec{u}}^* = \underset{\vec{u}}{\arg}\min \mathcal{E}(\vec{u}),\label{eq9}
\end{equation}
where ${\vec u}^{*}$ is the solution of the problem (\ref{eq1})--(\ref{eq2_2}), $\mathcal{E}$ is the total energy of system which is defined as
\begin{equation}
	\mathcal{E}\left( {\vec u} \right) = \frac{1}{2}\int_{\Omega} \left(\sigma_{ik}\varepsilon_{ik}\right) d V - \int_{\Omega} \left( {\vec f}\cdot {\vec u} \right) d V - \oint_{\partial \Omega} \left( {\vec T}\cdot {\vec u} \right) d S.   \label{eq10}
\end{equation}
Here $({\vec a}\cdot{\vec b})$ means dot product of two vector ${\vec a}$ and ${\vec b}$, ${\vec u}$ is constrained by boundary conditions (\ref{eq2_2}) for the displacements.

We will assume the displacements for the minimum of the total energy of the system coincide with displacements which are solutions of the problem described by Eqs.~(\ref{eq1})--(\ref{eq2_2}).

\section{Methods}

\subsection{Physical-Informed approach}
\subsubsection{Vanilla PINN}
According to the physics-informed neural network (PINN) approach~\cite{Raissi2019}, which stands on the shoulders of the universal approximation theorem~\cite{HORNIK1989359}, the PINN approximates the unknown solution ${\vec u}({\vec x})$ of the problem  (\ref{eq1})--(\ref{eq2_2}) by a deep neural network ${\vec u}^{{\bm \theta}}({\vec x})$, where ${\bm \theta}$ denote all trainable parameters of the network (weights and biases). Finding optimal parameters is an optimization problem, which requires the definition of a loss function such that its minimum gives the solution of the PDE. The physics-informed model is trained by minimizing the composite loss function which consists of the local residuals of the differential equations (\ref{eq1}) over the problem domain and its boundary as shown below:
\begin{equation}
	\mathcal{L}({\bm \theta}) = \lambda_{bc} \mathcal{L}_{bc}({\bm \theta}) + \mathcal{L}_{r}({\bm \theta}),\label{eq12}
\end{equation}

where
\begin{eqnarray}
	&& \mathcal{L}_{bc}\left(\bm \theta \right) = \frac{1}{N_{bc}} \sum_{i=1}^{N_{bc}} \sum_{l=1}^{N_l} \left| \mathcal{B}_l\left[{\vec u}^{{\bm \theta}} \right] \left({\vec x}_{i}^{(bc)}  \right)  \right|^{2}\label{eq6_2},\\
	&& \mathcal{L}_{r}\left(\bm \theta \right) = \frac{1}{N_r} \sum_{i=1}^{N_{r}}  \sum_{l=1}^{3} \left| \mathcal{R}_l\left[{\vec u}^{{\bm \theta}} \right] \left({\vec x}_{i}^{(r)} \right)\right|^{2},\label{eq7_2}\\
	&& \mathcal{B}_l[{\vec u}]:= {u}_l - {u}^{(bc)}_l \quad\text{or}\quad \mathcal{B}_l[{\vec u}]:= \sigma_{lk}n_{k} - T_{l}, \label{eq8_3} \\ 
	&& \mathcal{R}_l\left[{\vec u} \right]:= \frac{\partial \sigma_{lk}}{\partial x_k} + f_l.\label{eq8_2}
\end{eqnarray}
Here $\left\{{\vec x}_{i}^{(bc)}\right\}^{N_{bc}}_{i=1}$ and $\left\{{\vec x}_{i}^{(r)}\right\}^{N_{r}}_{i=1}$ are sets of points corresponding to boundary condition domain and PDE domain, respectively. These points can be the vertices of a fixed mesh or can be randomly sampled at each iteration of a gradient descent algorithm. $\mathcal{B}_l$ is a boundary operator corresponding to $l$th boundary conditions at the given boundary point, $\mathcal{R}_l$ is a residual operator which is determined by the $l$th equation of PDEs (\ref{eq1}).

All required gradients, derivations (for example, $\sigma_{lk}$ and $\partial\sigma_{lk} / \partial x_k$) w.r.t. input variables ($\vec x$) and network parameters $\theta$ can be efficiently computed via automatic differentiation~\cite{Griewank2008} with algorithmic accuracy, which is defined by the accuracy of computation system. The hyperparameter $\lambda_{bc}$ allows for separate tuning of the learning rate for each of the loss terms in order to improve the convergence of the model~\cite{Wang2021,WANG2022110768}.\\

The optimization problem can be defined as follows
\begin{equation}
	{\bm \theta}^* = \underset{{\bm \theta}}{\arg}\min \mathcal{L}({\bm \theta}),\label{eq17}
\end{equation}
where ${\bm \theta}^*$ are optimal parameters of the neural network which minimize the discrepancy between the exact unknown solution ${\vec u}$ and the approximate one ${\vec u}_{\pmb \theta^*}$.

\subsubsection{Energy PINN}

Following the deep energy method (DEM)~\cite{Samaniego2020} the energy of deformation of the mechanical system (\ref{eq10}) can be used as the energy term of the loss function (\ref{eq12}) instead of residual loss term $\mathcal{L}_r(\bm \theta)$. Then, the following form for the loss function is obtained
\begin{equation}
	\mathcal{L}({\bm \theta}) = \lambda_{bc} \mathcal{L}_{bc}({\bm \theta}) + \mathcal{L}_{en}({\bm \theta}),\label{eq18}
\end{equation}
where
\begin{align}
	\mathcal{L}_{en}({\bm \theta}) & = \frac{1}{2}\sum_{i=1}^{N_{en}} w_i \left(\sigma_{lk}\left({\vec x}_{i}^{(en)}\right)\varepsilon_{lk}\left({\vec x}_{i}^{(en)}\right)\right) \Delta V_i\notag\\
	& - \sum_{i=1}^{N_{en}} w_i \left( {\vec f}\cdot {{\vec u}^{{\bm \theta}}(x_i)} \right) \Delta V_i - \sum_{i=1}^{N_{bc}} w_i \left( {\vec T}\cdot {\vec u}_{\bm \theta}\left({\vec x}_{i}^{(bc)}\right) \right) \Delta S_i.\label{eq19}
\end{align}
Here $\left\{{\vec x}_{i}^{(en)}\right\}^{N_{en}}_{i=1}$ are set of points corresponding to the domain of body are taken for the integration, $\Delta V_i$ is the volume of $i$th cell of the body, $\Delta S_i$ is square of $i$th cell of the surface of the body, $w_i$ is the weight coefficient of given integration scheme in $i$th point (in our calculations we taken the Simpson's rule of integration).

Note, in this case of the loss function, the term $\mathcal{L}_{bc}$ is generated by the elements corresponding to boundary conditions (\ref{eq2_2}), but is not for (\ref{eq2_1}) because last conditions are involved in the loss term (\ref{eq19}) as an energy of external forces.\\
The optimization problem is still given by the Eq.~(\ref{eq17}).

\subsection{Artificial Neural Network Architectures}

A generic deep neural network with $L$ layers, can be expressed by the composition of $L$ functions $f_i({\vec z}_i; {\bm \theta}_i)$~\cite{Caterini2018,Cuomo2022}
(where ${\vec z}_i$ are the state variables, and ${\bm \theta}_i$ is the set of parameters of the $i$th layer). The function {of interest} ${\vec u}_{\bm \theta}$ can be {represented} in the following form
\begin{equation}
	{\vec u}_{\bm \theta}({\vec x}) = f_{L}\circ f_{L-1}\circ \dots \circ f_1({\vec x}),\label{eq20}
\end{equation}
where functions $f_{i}$ {are determined with inner product of two spaces} $E_i$ and $H_i$, so that $f_i \in E_i\times H_i$; $\circ$ is operator of  function composition which reads for function $f_1$ and $f_2$ as $f_2\circ f_1({\vec x}) = f_2\left(f_1({\vec x})\right) $.

\subsubsection{Multi-layer perceptron}

As in the original vanilla PINN~\cite{Raissi2017,Raissi2017_2}, the feed-forward neural networks (multi-layer perceptron (MLP)) play the key role in the representation of the solutions for the problems in our paper. This neural network is a set of neurons collected in layers which execute calculations sequentially, layer by layer. Data flow through the MLP in a forward direction without loopback. The neurons of feed-forward neural networks in neighboring layers are connected, whereas neurons of a single layer are not directly related. Every neuron is a collection of sequence mathematical operations which are the sum of weighted neuron's inputs and a bias factor, {and appliance an activation function to the result}. The action of one layer of the MLP can be written as  
\begin{equation}
	f_i({\vec z}_i; {\vec W}_i, {\vec b}_i) = \alpha_i \left({\vec W}_i \cdot {\vec z}_i + {\vec b}_i\right), \label{eq21}
\end{equation}
where $z_i$ is output of th $(i-1)$th layer, ${\vec W}_i$ and ${\vec b}_i$ are the the weight and the bias of $i$th layer, respectively, $\alpha_i$ is a non-linear activation function. The general representation (\ref{eq20}) of artificial neural network for the MLP with identical activation function for all layers can also be rewritten in conformity with the notation~\cite{Mishra2021}

\begin{equation}
	{\vec u}_{\bm \theta}({\vec x}) = C_{L}\circ \alpha \circ C_{L-1}\circ \dots \circ \alpha \circ C_1({\vec x}),\label{eq22}
\end{equation}
where for any $k\in [1,2,\cdot, L]$ the $C_k$ is defined with
\begin{equation}
	C_k({\vec x}_k) = {\vec W}_k {\vec x}_k +  {\vec b}_k.\label{eq23}
\end{equation}
Thus this MLP consists of input and output layers, and $(L-2)$ hidden layers.

\begin{figure}[t!]
	\centering
	\includegraphics[width=1\textwidth]{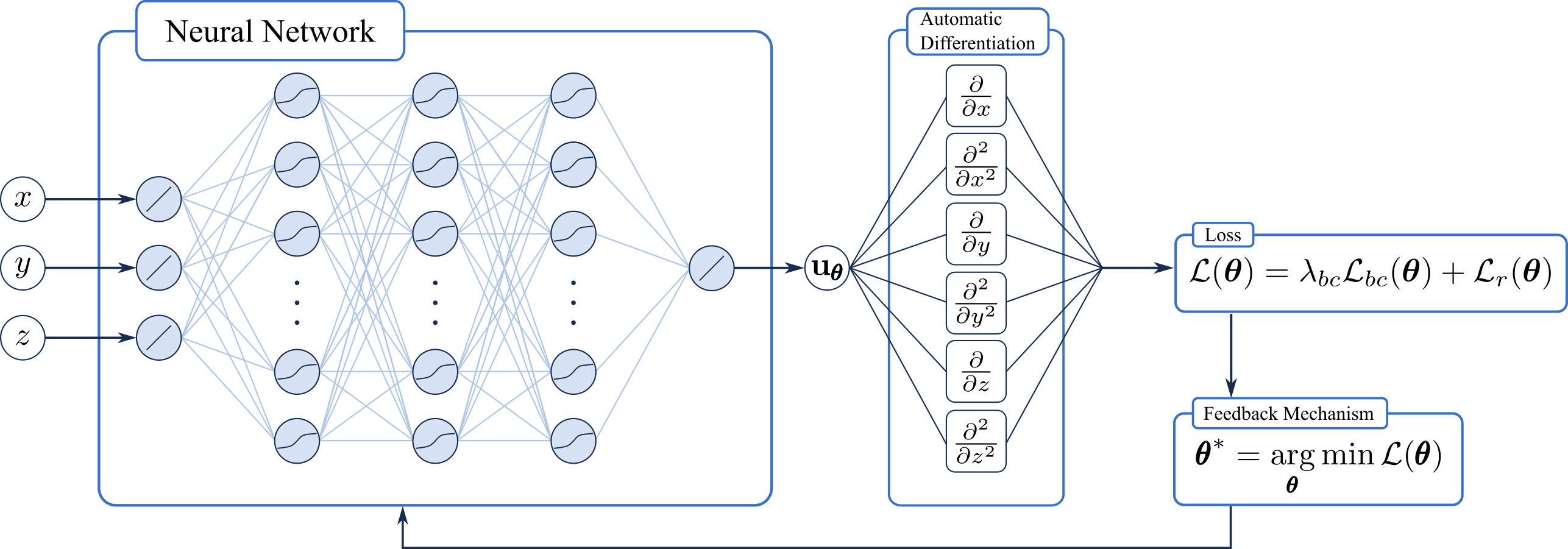}
	\caption{Schematic diagram of the general PINN method.}
	\label{fig1}
\end{figure}

\subsubsection{Separable PINN}

In most part of our experiments, we used a novel network architecture, called separable PINN (SPINN), which facilitates significant forward-mode automatic differentiation (AD) for more efficient computation~\cite{cho2023}. The brilliant idea behind this approach is using the classic separation of variables method (the Fourier method) for the special representation of the neural network. SPINN uses a per-axis basis instead of point-wise processing in vanilla PINNs, which decreases the number of network forward passes. In addition, training time, evaluating time and memory costs are significantly reduced for the SPINN in comparison to PINN, for which these parameters are grown exponentially along with the grid resolution increasing.\\

SPINN consists of $d$ body-networks (which are usually presented by MLPs), each of which takes an individual one-dimensional coordinate component as an input from $d$-dimensional space. Each body-network ${\vec f}_{{\bm \theta}_i}:\mathbb{R}\rightarrow \mathbb{R}^{r}$ (parameterized
by ${{\bm \theta}_i}$) is a vector-valued function which transforms the coordinates of $i$th axis into a $r$-dimensional feature representation (functional basis, or set of eigenfunctions).

The predictions of SPINN are computed by basis functions merging which is defined by following
\begin{equation}
	{u}^{{\bm \theta}}({\vec x}) = \sum\limits_{j=1}^{r} \prod\limits_{i=1}^{d} {f}^{{\bm \theta}_i}_{j} (x_i),\label{eq24}
\end{equation}
where ${u}^{{\bm \theta}}:\mathbb{R}^{d}\rightarrow \mathbb{R}$ is the predicted solution function, $x_i \in \mathbb{R}$ is coordinate of $i$th axis, and ${f}^{{\bm \theta}_i}_{j}$ is the $j$th element of basis ${\vec f}^{{\bm \theta}_i}$ $\left({\vec f}^{{\bm \theta}_i} = \left({f}^{{\bm \theta}_i}_{1},{f}^{{\bm \theta}_i}_{2},\dots,{f}^{{\bm \theta}_i}_{r}\right)\right)$.

\begin{figure}[t!]
	\centering
	\includegraphics[width=1\textwidth]{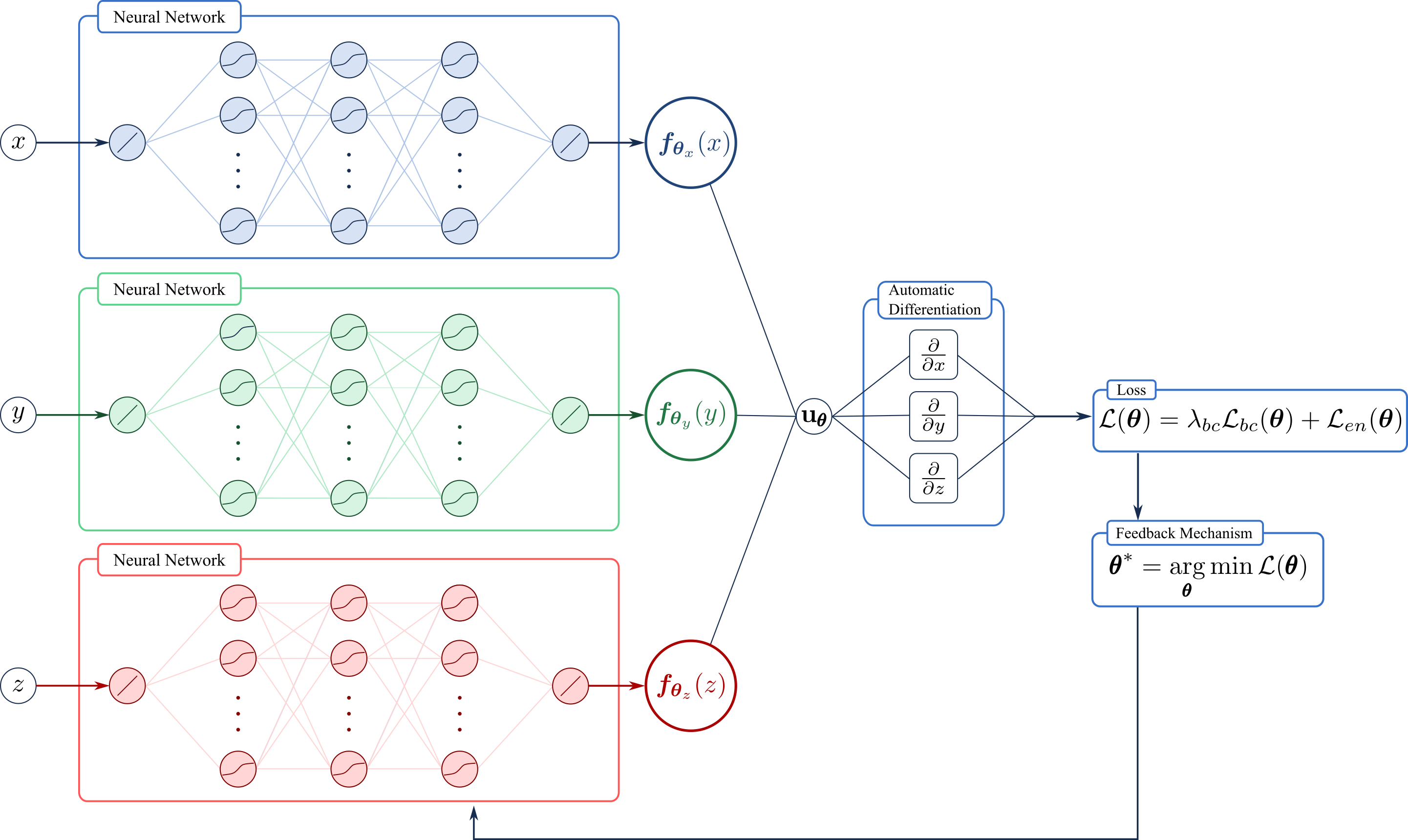}
	\caption{Schematic diagram of the deep energy method for the SPINN approach.}
	\label{fig2}
\end{figure}

A schematic diagram of the general PINN method is shown in Figure~\ref{fig1}.

\section{Numerical experiments}\label{experiments}

To evaluate the accuracy of the approximate solution obtained with the help of the PINN method, the values of the solution of the linear elasticity problem predicted by the neural networks at given points are compared with the values calculated on the basis of the classical high-precision numerical method (FEM). The relative total ${\mathbb L}_2$ error of prediction is taken as a measure of accuracy, which can be expressed with the following relation
\begin{equation}\label{eq61}
	\epsilon_{\rm error} = \left\{\frac{1}{N_e} \sum^{N_e}_{i=1} \left[{\vec u}^{\pmb \theta}({\vec x}_i) - {\vec u}({\vec x}_i)\right]^2 \right\}^{1/2} \times \left\{\frac{1}{N_e} \sum^{N_e}_{i=1} \left[{\vec u}({\vec x}_i)\right]^2 \right\}^{-1/2},
\end{equation}
where $\left\{{\vec x}_{i}\right\}^{N_{e}}_{i=1}$ is the set of evaluation points taken from the domain $\Omega$ of body, ${\vec u}_{\pmb \theta}$ and ${\vec u}$ are the predicted and reference solutions respectively.

Relative total ${\mathbb L}_2$ errors were analyzed for the {displacement components} separately, for the magnitude of the displacements, which are given by 
\begin{equation}
	u_{\rm m} = \left(u_x^2 + u_y^2 + u_z^2\right)^{1/2},\label{eq26}
\end{equation}
and von Mises stress
\begin{equation}
	\sigma_{\rm vM} =  \left[\frac{3}{2} \sigma_{ij}\sigma_{ij} -\frac{1}{2} \left(\sigma_{kk}\right)^2\right]^{1/2}.
\end{equation}

For our experiments we used Pytorch~\cite{Paszke2019} version 1.12.1 {as backend for} Modulus, and Jax, which was used with code based on the source code is available at \url{https://github.com/stnamjef/SPINN} of paper~\cite{cho2023}. The training was carried out on a node with Nvidia Tesla V100 GPU.

Reference values were obtained using OpenFOAM at CPU. There are presented the {timespan} of calculations of reference values in tables~2 and 5 {considering} the mesh generation time.

\begin{table}[h!]
\centering
\label{table1}
\begin{tabular}{lcc}
	\toprule
	\head{3cm}{Description} \rule[-1ex]{0pt}{3.5ex}  & \head{3cm}{Values for the beam problem} & \head{3cm}{ Values for the thin-walled angle problem}\\
	\midrule
	{$E$ (Young modulus) (Pa)} & \multicolumn{2}{c}{$2.1 \times10^{11}$} \\		
	\midrule
	{$\nu$ (Poisson ratio)} & \multicolumn{2}{c}{0.3} \\
	\midrule
	\head{6cm}{$T$ (traction on external boundary) (Pa)} & $10^4$ & $2.5 \times 10^4$ \\
	\midrule
	$g$ (gravitational acceleration) (m/s$^2$) & \multicolumn{2}{c}{9.81}\\	
	\midrule
	$\rho$ (density of body material) (kg/m$^3$) & \multicolumn{2}{c}{$7.8 \times 10^3$}\\
	\bottomrule
\end{tabular}
\caption{Parameters for the problems {(material is steel)}}
\end{table}

For improving the training convergence and accuracy~\cite{Modulus_Linear_Elasticity,wang2023experts} we transformed the original Eqs.(\ref{eq1})--(\ref{eq5}) to the non--dimensionalized forms. The non--dimensionalized variables are defined as follows
\begin{align}
	& \tilde{x}_i = \frac{x_i}{L_c},\quad \tilde{u}_i = \frac{u_i}{U_c},\quad \tilde{\lambda} = \frac{\lambda}{\mu_c},\quad \tilde{\nu} = \frac{\mu}{\mu_c}, \label{eq28}
\end{align}
where $L_c$ is the characteristic length, $U_c$ is the characteristic displacement, and $\mu_c$ is the non--dimensionalizing shear modulus. The non--dimensionalized state of a static equilibrium problem can be written as follows:
\begin{align}\label{eq29}
	&\frac{\partial \tilde{\sigma}_{ik}}{\partial \tilde{x}_k} + \tilde{f}_i = 0,\,\quad \tilde{{\vec x}}\in \tilde{\Omega}\\
	&\tilde{\sigma}_{ik}n_{k} = \tilde{T}_{i}, \quad \tilde{\vec x} \in \partial \tilde{\Omega},\label{eq30}\\
	&\tilde{\vec u} = \tilde{\vec u}^{(b)}, \quad \tilde{\vec x} \in \partial \tilde{\Omega},\label{eq31}
\end{align}

where the non--dimensionalized forces and stress tensor are
\begin{align}
	& \tilde{\sigma}_{ik} = \frac{L_c}{\mu_c U_c} {\sigma}_{ik},\quad \tilde{f}_{i} = \frac{L_c^2}{\mu_c U_c} {f}_{i}, \quad \tilde{T}_i = \frac{L_c}{\mu_c U_c} T_i.  \label{eq32}
\end{align}
The values $\tilde{\vec u}^{(b)}$, $\tilde\Omega$ and $\partial \tilde\Omega$ are ${\vec u}^{(b)}$, $\Omega$ and $\partial \Omega$ were taken in transformed coordinates.
Non--dimensionalized stress--displacement and strain--displacement relations are obtained by multiplying of equations (\ref{eq3}) and (\ref{eq4}) by $L_c / (\mu_c U_c)$ and $L_c / U_c$, respectively, and have the following forms
\begin{align}
	& \tilde{\sigma}_{ik} = \tilde{\lambda} \tilde{\varepsilon}_{ll} \delta_{ik} + \tilde{\mu} \tilde{\varepsilon}_{ik}, \label{eq33}\\
	& \tilde{\varepsilon}_{ik} = \frac{1}{2} \left( \frac{\partial \tilde{u}_i}{\partial \tilde{x}_{k}} + \frac{\partial \tilde{u}_{k}}{\partial \tilde{x}_{i}}\right).\label{eq34}
\end{align}

\subsection{Beam}

\begin{figure}[ht!]
	\centering
	\includegraphics[width=0.5\textwidth]{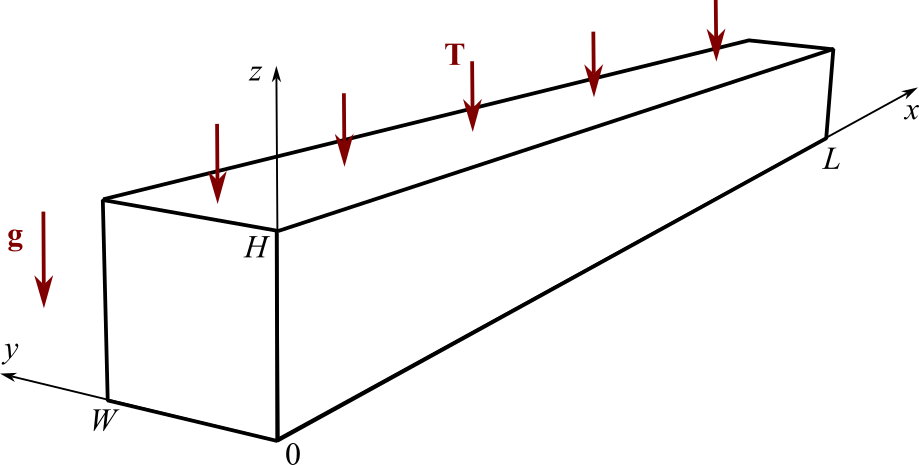}
	\caption{Beam with uniform load $\vec T$ in the gravity force $\vec g$}
	\label{fig3}
\end{figure}

Consider a beam made of an isotropic homogeneous elastic material with length $L = 1$ m, width $W = 0.1$ m and height $H = 0.1$ m (see Fig.~\ref{fig3}). So the coordinates of points of the beam ($x,y,z$) is satisfied to relations $x\in [0,1]$, $y\in [0,0.1]$, $z\in [0,0.1]$). The beam is fixed on boundary $x=0$, and is loaded by pressure ${\vec T} = - T {\vec z}_0$ on boundary $z=0.1$ m and gravity force ${\vec f} = - \rho g {\vec z}_0 $ in each point of the body. Parameters of materials and magnitudes of forces are given in Table~1.

The undimensionalization parameters were chosen as $L_c = L$, $\mu_c = 0.01 \mu$ and $U_c = 10^{-4}$ m.

\begin{align}
	& \frac{\partial \tilde{\sigma}_{xx}}{\partial \tilde{x}} + \frac{\partial \tilde{\sigma}_{xy}}{\partial \tilde{y}} + \frac{\partial \tilde{\sigma}_{xz}}{\partial \tilde{z}} = 0, \notag\\
	& \frac{\partial \tilde{\sigma}_{xy}}{\partial \tilde{x}} + \frac{\partial \tilde{\sigma}_{yy}}{\partial \tilde{y}} + \frac{\partial \tilde{\sigma}_{yz}}{\partial \tilde{z}} = 0, \notag\\
	& \frac{\partial \tilde{\sigma}_{xz}}{\partial \tilde{x}} + \frac{\partial \tilde{\sigma}_{yz}}{\partial \tilde{y}} + \frac{\partial \tilde{\sigma}_{zz}}{\partial \tilde{z}} - \tilde{\rho} g= 0,  \notag\\
	& \tilde{u}_x = 0,\quad \tilde{u}_y = 0,\quad \tilde{u}_z = 0 \quad \text{in} \quad \tilde{x}=0, \notag\\
	& \tilde{\sigma}_{zz} = -\tilde{T}, \quad \tilde{\sigma}_{xz} = 0, \quad \tilde{\sigma}_{yz} = 0 \quad \text{in} \quad \tilde{z}=0.1, \notag\\
	&\tilde{\sigma}_{ik} = 0,  \quad \text{in} \quad \tilde{z}=0, \quad  \tilde{x}=1\quad  \tilde{y}=0,\quad  \tilde{y}=0.1.
\end{align}
Here $\tilde{\rho} = {L_c^2 \rho}/{(\mu_c U_c)}$ and  $\tilde{T} = {L_c T}/{(\mu_c U_c)}$.

\subsubsection{PINN}
The PINN approach was {implemented on the basis of Modulus from Nvidia} with recommended parameters for the linear elasticity problems~\cite{Modulus_Linear_Elasticity}. {Each value of displacements and stress components are predicted with the help of} the individual feedforward neural network (9 MLPs in the sum). Each MLP has 6 hidden layers with 512 neurons per layer. The activation function of these MLPs is ``swish'' function $[\text{swish}(x) = x\ \text{sigmoid}(x)]$.

The networks had been trained via stochastic gradient descent using the Adam optimizer~\cite{Kingma2014}. The initial value of the learning rate is $10^{-3}$. Training consisted of 2,000,000 epochs for the method based on PDE. Weight factors $\lambda_{bc}$ in the loss function in Eq.~(\ref{eq12}) is set to 10. The training is based on 50,000 collocation points for the inner area of the beam (volumetric points) and 1664 points for the surface  (surface points).

\subsubsection{SPINN in the case of the loss based on PDE}
The values of displacements and stresses are approximated {with} two different neural networks. First of which predicts three components of vector $\tilde{\vec u}$, second neural network predicts six components of tensor $\tilde{\sigma}_{ik}$. Due to this circumstance the loss function (\ref{eq12}) must be added by term
\begin{equation}
	\mathcal{L}_{u\sigma}({\bm \theta}) = \frac{1}{N_r} \sum_{j=1}^{N_{r}} \sum_{i=1}^{3}  \sum_{k=1}^{3} \left| \tilde{\sigma}_{ik}^{{\bm \theta}} - \left\{\tilde{\lambda} \frac{\partial \tilde{u}^{{\bm \theta}}_{l}}{\partial \tilde{x}_{l}} \delta_{ik} +  \frac{\tilde{\mu}}{2} \left( \frac{\partial \tilde{u}^{{\bm \theta}}_i}{\partial \tilde{x}_{k}} + \frac{\partial \tilde{u}^{{\bm \theta}}_{k}}{\partial \tilde{x}_{i}}\right)\right\} \right|^{2},\label{eq36}
\end{equation}
and $\sigma_{ik}$ {should} be changed by $\tilde{\sigma}_{ik}^{{\bm \theta}}$ in all expressions (\ref{eq6_2})--(\ref{eq8_2}). We will name this approach SPINN--PDE in our paper for brevity. \\

\subsubsection{SPINN conjunction with DEM}

In the SPINN approach, {each} neural network (for displacements or stresses predictions) consists of the three body--networks. Each body--network is a single MLP which has 3 hidden layers with 64 neurons per layer and 64 output values that play roles of the basis functions. The activation function of these MLPs is ``swish'' function. We will name this approach SPINN--DEM in our paper for brevity. \\

The networks had been trained via stochastic gradient descent using the Adam optimizer~\cite{Kingma2014} and scheduler, which reduces the learning rate by 5\% every 5,000 epochs (decay--rate is 0.95). The initial value of the learning rate is $10^{-3}$. Training consisted of 500,000 epochs for the method based on PDE and 100,000 epochs for the DEM. Weight factors $\lambda_{bc}$ in the loss function in Eq.~(\ref{eq12}) is set to 100. The training is based on $N_r=32^3$ collocation points ($32$ points for every axis) which are resampled every 100 epochs, except the SPINN--DEM approach, where the points were fixed. Evaluating the results of training carried out every 1000 epochs. Every training for this approach is performed 7 times with different seeds.

\subsubsection{Comparison of the results}
The experiments were carried out for all the above-mentioned methods: PINN based on PDE, SPINN based on PDE, and SPINN based on the DEM. In table~{2} we present the {timespan} of one experiment, speed of training, and GPU memory consumption for these approaches. It can be seen that SPINN bypasses the PINN by a large margin both in speed and in memory using efficiency. In table~{3} we also report the accuracy for the beam problem achieved using different methods for $u_{x,y,z}$, $u_{\rm m}$ and $\sigma_{\rm vM}$. As can be seen from this table the best result is achieved with SPINN based on DEM for all values while PINN and SPINN demonstrate close results. 

In the table~{4} are given times to achieve ``engineering'' accuracy  ${\mathbb T}_{5\%}$ for the different values obtained by different approaches (in seconds). {We defined “engineering” accuracy as} relative error which should be less or equal 5$\%$ (${\mathbb L}_2 \leq 0.05$). This table demonstrates the significant superiority of the SPINN based on DEM over the other methods in the speed of achieving the target accuracy of values.

\begin{figure}[ht!]
	\centering
	\includegraphics[width=1\textwidth]{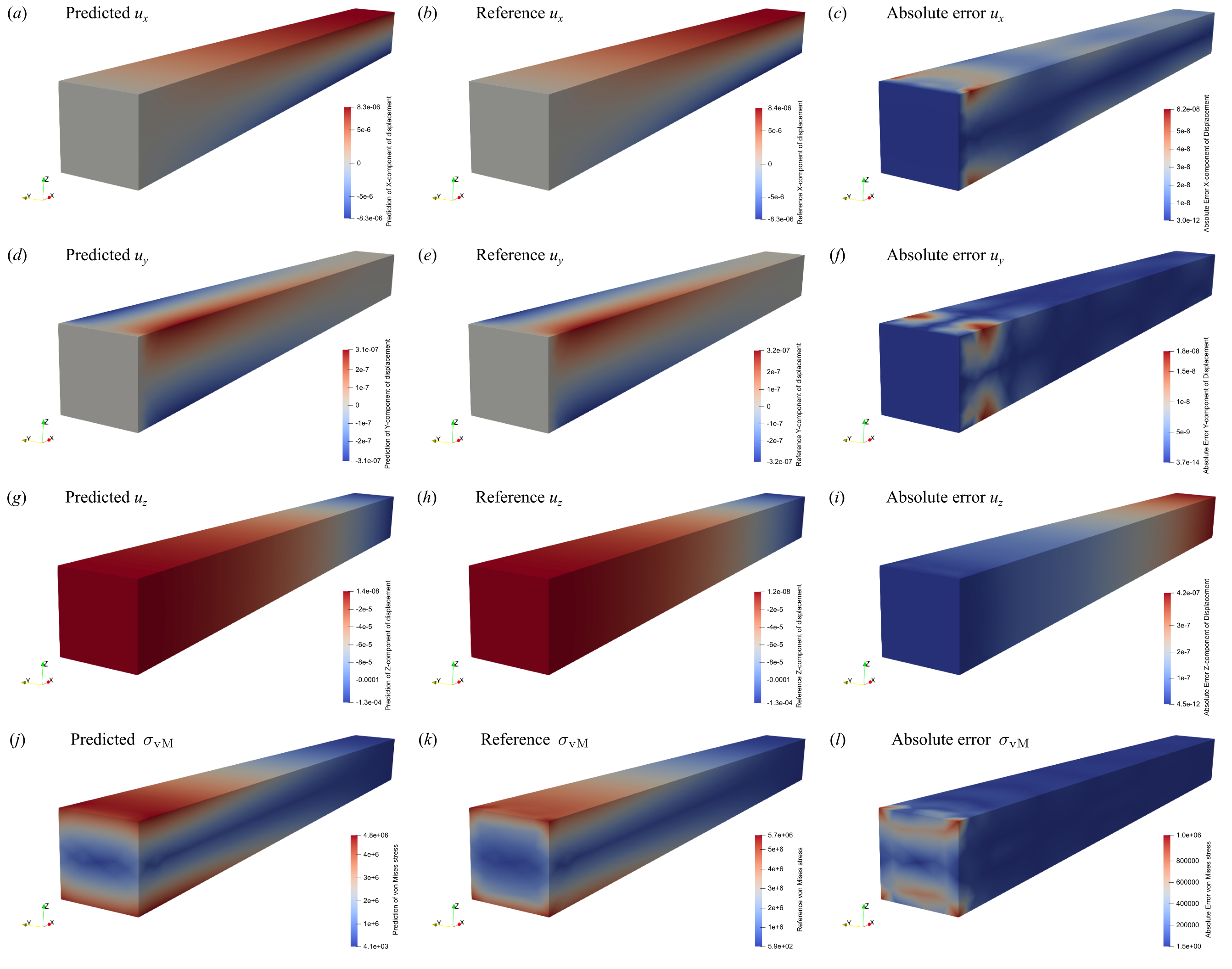}
	\caption{Beam problem: (a), (d), (g) and (j) are the predictions of a trained physics-informed neural network; (b), (e), (h) and (k) are reference solutions; (c), (f), (i) and (l) are absolute values of differences between reference and predicted solutions, for the value $u_x$, $u_y$, $u_z$ and $\sigma_{\rm vM}$, respectively.}
	\label{fig4a}
\end{figure}

The displacements and value of von Mises stress predicted by SPINN based on DEM, reference values obtained by FEM, and absolute errors are presented in Fig.~\ref{fig4a}. As you can see, SPINN gives distributions of values very close to the reverences. The most pronounced difference is observed in the von Mises stress at the place of fixing of the beam. However, distributions found by SPINN is more smooth in the points of fixation than distributions given by FEM. In our opinion, this smoothness of solutions that are obtained by SPINN indicates that this solution is closer to the real physical distributions than the one given by the FEM.

\begin{table}
	\centering
	\begin{tabular}{lccc}
		\toprule
		{Method} \rule[-1ex]{0pt}{3.5ex}  & Runtime (s) & Speed (epochs/s) & GPU Memory (MB)\\
		\midrule
		{FEM (based on PDE (\ref{eq6})--(\ref{eq8}), at CPU) } & 4  & --- & --- \\
		\midrule
		{PINN (Modulus, loss (\ref{eq12})) } & 277515 & 7 & 14501 \\		
		\midrule
		{SPINN--PDE} & $1869.86\pm3.47$ & 283 & 733 \\
		\midrule
		{SPINN--DEM} & $356.54\pm3.70$
		& 295 & 723 \\
		\bottomrule
	\end{tabular}\label{table2}
\caption{Beam problem: {Timespan} of one experiment, speed of training and GPU memory consumption for different approach}
\end{table}

\begin{table}
\centering
\begin{tabular}{lccccc}
	\toprule
	{Method} \rule[-1ex]{0pt}{3.5ex}  & {${\mathbb L}_2[u_x]$} & {${\mathbb L}_2[u_y]$} & {${\mathbb L}_2[u_z]$} & {${\mathbb L}_2[u_{\rm m}]$} & {${\mathbb L}_2[\sigma_{\rm vM}]$}\\
	\midrule
	{PINN (Modulus, loss (\ref{eq12})) } & 0.0327 & 0.4465 & 0.0348 & 0.0348 & --- \\		
	\midrule
	{SPINN--PDE} & $\underset{\pm0.0034}{0.0155}$
	& $\underset{\pm0.0123}{0.3252}$ 	&$\underset{\pm0.0034}{0.0129}$  
	&$\underset{\pm0.0034}{0.0129}$  
	&$\underset{\pm0.0005}{0.0719}$  
	\\
	\midrule
	{SPINN--DEM} & $\underset{\pm0.0010}{0.0025}$
	& $\underset{\pm0.0089}{0.0281}$	& $\underset{\pm0.0008}{0.0030}$  
	& $\underset{\pm0.0008}{0.0030}$  
	& $\underset{\pm0.0049}{0.0419}$ \\
	\bottomrule
\end{tabular}\label{table4}
\caption{Beam problem: Relative ${\mathbb L}_2$ errors obtained by different approaches for the different values}
\end{table}

\begin{table}
\centering
\begin{tabular}{lccccc}
	\toprule
	{Method} \rule[-1ex]{0pt}{3.5ex}  & {${\mathbb T}_{5\%}[u_x]$} & {${\mathbb T}_{5\%}[u_y]$} & {${\mathbb T}_{5\%}[u_z]$} & {${\mathbb T}_{5\%}[u_{\rm m}]$} & {${\mathbb T}_{5\%}[\sigma_{\rm vM}]$}\\
%	\midrule
%	{PINN (Modulus, loss (\ref{eq12})) } & & & & & \\		
	\midrule
	{SPINN--PDE} & $\underset{\pm57.44}{441.48}$ & --- & $\underset{\pm48.77}{426.84}$ & $\underset{\pm48.77}{426.84}$ & --- \\
	\midrule
	{SPINN--DEM} & $\underset{\pm2.67}{9.75}$ & $\underset{\pm19.84}{55.08}$ & $\underset{\pm2.67}{9.75}$ & $\underset{\pm2.67}{9.75}$ & $\underset{\pm20.97}{112.54}$ \\
	\bottomrule
\end{tabular}\label{table3}
\caption{Beam problem: Times to achieve ``engineering'' accuracy  ${\mathbb T}_{5\%}$ for the different values obtained by different approaches (in seconds). ``Engineering'' accuracy is defined by relative error which should be less or equal 5$\%$ (${\mathbb L}_2 \leq 0.05$).}
\end{table}

%\begin{figure}[t!]
%	\centering
%	\includegraphics[width=1\textwidth]{figs/fig4_b.png}
%	\caption{(a) is prediction of a trained physics-informed neural network, (b) is reference solution, (c) are absolute values of difference between reference and predicted solutions, for the value $u_y$.}
%	\label{fig4b}
%\end{figure}
%
%\begin{figure}[t!]
%	\centering
%	\includegraphics[width=1\textwidth]{figs/fig4_c.png}
%	\caption{(a) is prediction of a trained physics-informed neural network, (b) is reference solution, (c) are absolute values of difference between reference and predicted solutions, for the value $u_z$.}
%	\label{fig4c}
%\end{figure}
%
%\begin{figure}[t!]
%	\centering
%	\includegraphics[width=1\textwidth]{figs/fig4_d.png}
%	\caption{(a) is prediction of a trained physics-informed neural network, (b) is reference solution, (c) are absolute values of difference between reference and predicted solutions, for the value $\sigma_{\rm vM}$.}
%	\label{fig4d}
%\end{figure}

\subsection{Thin-walled angle}

\begin{figure}[t!]
	\centering
	\includegraphics[width=0.5\textwidth]{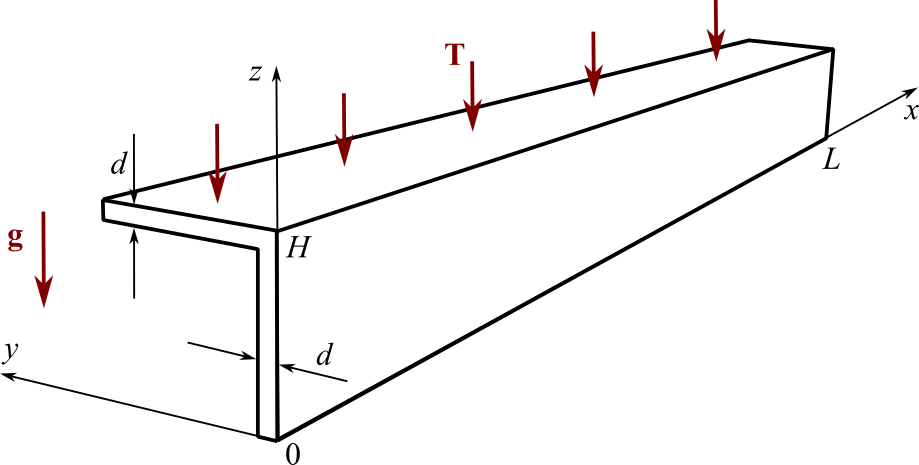}
	\caption{Thin-walled angle with uniform load $\vec T$ in the gravity force $\vec g$}
	\label{fig3_1}
\end{figure}

Consider a thin-walled angle made of an isotropic homogeneous elastic material with length $L = 1$ m, width $W = 0.06$ m, height $H = 0.06$ m, and the wall thickness $d=6\times 10^{-3}$ m (see Fig.~\ref{fig3_1}). So the coordinates of points of the angle ($x,y,z$) is satisfied to relations
\begin{equation}
	x\in [0,1],\quad (y, z) \in [0,0.06]\times[0,0.06],\quad (y, z) \notin [0,0.054]\times[0,0.054]. \notag
\end{equation}
The angle is fixed on boundary $x=0$, and is loaded by pressure ${\vec T} = - T {\vec z}_0$ on boundary $z=0.06$ m and gravity force ${\vec f} = - \rho g {\vec z}_0 $ in each point of the body. Parameters of materials and magnitudes of forces are given in Table~1.\\

The {non-dimensional} parameters were chosen as $L_c = L$, $\mu_c = \mu$ and $U_c = 10^{-4}$ m.

Note, that no one method (vanilla PINN and SPINN) based on the partial differential equations (\ref{eq12}) {did not converge} to reference solutions in our experiments.

In the SPINN--DEM approach every element of the displacements vector is predicted by an individual neural network. Each such neural network consists of the three body--networks. Each body--network is a single MLP which has 5 hidden layers with 64 neurons per layer and 64 output values which play roles of the basis functions. The activation function of these MLPs is ``swish'' function.

The networks had been trained via stochastic gradient descent using the Adam optimizer~\cite{Kingma2014} and scheduler, which reduces the learning rate by 5\% every 5,000 epochs (decay--rate is 0.95). The initial value of the learning rate is $10^{-3}$. Training consisted of 200,000 epochs. Weight factors $\lambda_{bc}$ in the loss function in Eq.~(\ref{eq18}) is set to 1000. The domain of integration was separated on two domains: $[0,1]\times[0,0.06]\times[0,0.54]$ and $[0,1]\times[0,0.6]\times[0.054,0.06]$. The collocation points were taken on uniform grids. The number of points for the first and second domains is $513\times9\times65$ and  $513\times65\times9$, respectively. Evaluating the results of training carried out every 1000 epochs. Every training for this approach is performed 7 times with different seeds.

In table~{5} we present the duration of one experiment, speed of training, and GPU memory consumption for training of SPINN. It can be seen that SPINN has speed and memory efficiency close to the SPINN in the previous problem. In table~{6} we also report the accuracy for the problem achieved using this method for $u_{x,y,z}$, $u_{\rm m}$ and $\sigma_{\rm vM}$. As can be seen from this table the results achieved with SPINN based on DEM are relatively good. 

In the table~{7} are given times to achieve ``engineering'' accuracy  ${\mathbb T}_{5\%}$ for the different values obtained by this approach (in seconds). This table demonstrates the relatively good convergence time of SPINN to target accuracy.

\begin{figure}[t!]
	\centering
	\includegraphics[width=1\textwidth]{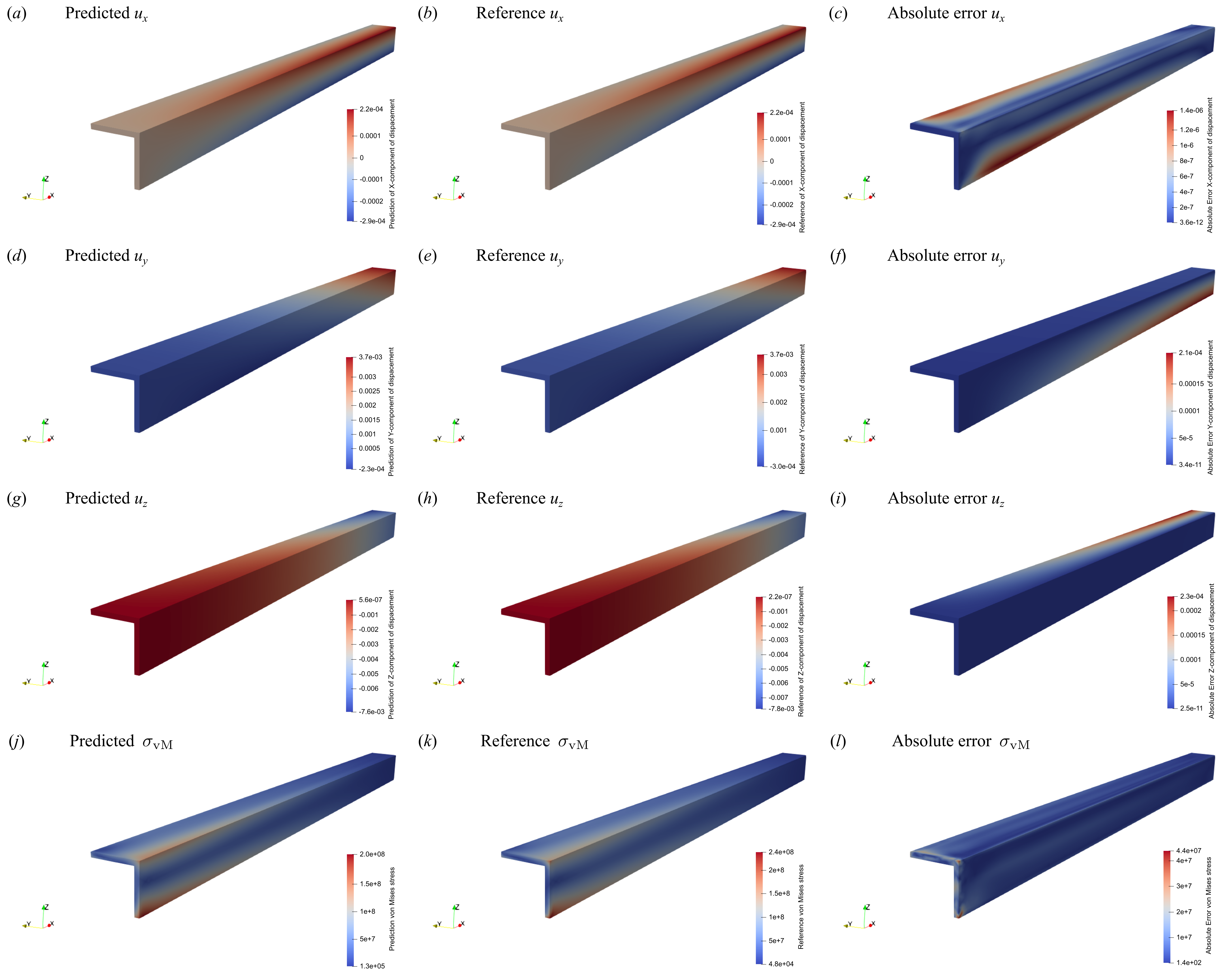}
	\caption{Thin-walled angle problem: (a), (d), (g) and (j) are prediction of a trained physics-informed neural network; (b), (e), (h) and (k) are reference solutions; (c), (f), (i) and (l) are absolute values of differences between reference and predicted solutions, for the value $u_x$, $u_y$, $u_z$ and $\sigma_{\rm vM}$, respectively..}
	\label{fig5a}
\end{figure}

The displacements and value of von Mises stress predicted by SPINN based on DEM, reference values obtained by FEM, and absolute errors are presented in Fig.~\ref{fig5a}. As before shown for the beam problem, SPINN gives distributions of values very close to the reverences once, and the most pronounced differences are laid in the von Mises stress at the place of fixing of the beam. The distributions found by SPINN look like more smooth in the {fixed points} than distributions given by FEM, and this smoothness of solutions indicates that this solution is closer to the real physical distributions than the one given by the FEM.

\begin{table}
\centering
\begin{tabular}{lccc}
	\toprule
	{Method} \rule[-1ex]{0pt}{3.5ex}  & Runtime (s) & Speed (epochs/s) & GPU Memory (MB)\\
	\midrule
	{FEM (based on PDE (\ref{eq6})--(\ref{eq8}), at CPU) } & 35 & --- & --- \\
	\midrule
	{SPINN--DEM} & $2159.64\pm24.1
	$  & 95 & 989 \\
	\bottomrule
\end{tabular}\label{table5}
\caption{Thin-walled angle problem: Duration of one experiment, speed of training and GPU memory consumption for different approach}
\end{table}

\begin{table}[t!]
\centering
\begin{tabular}{lccccc}
	\toprule
	{Method} \rule[-1ex]{0pt}{3.5ex}  & {${\mathbb L}_2[u_x]$} & {${\mathbb L}_2[u_y]$} & {${\mathbb L}_2[u_z]$} & {${\mathbb L}_2[u_{\rm m}]$} & {${\mathbb L}_2[\sigma_{\rm vM}]$}\\
	\midrule
	{SPINN--DEM} & $\underset{\pm0.0003}{0.0038}$
	& $\underset{\pm0.0037}{0.0389}$	& $\underset{\pm0.0017}{0.0205}$  
	& $\underset{\pm0.0014}{0.0173}$  
	& $\underset{\pm0.0034}{0.0807}$ \\
	\bottomrule
\end{tabular}\label{table6}
\caption{Thin-walled angle problem: Relative ${\mathbb L}_2$ errors obtained by SPINN--DEM approach}
\end{table}

\begin{table}[t!]
\centering
	\begin{tabular}{lccccc}
		\toprule
		{Method} \rule[-1ex]{0pt}{3.5ex}  & {${\mathbb T}_{5\%}[u_x]$} & {${\mathbb T}_{5\%}[u_y]$} & {${\mathbb T}_{5\%}[u_z]$} & {${\mathbb T}_{5\%}[u_{\rm m}]$} & {${\mathbb T}_{5\%}[\sigma_{\rm vM}]$}\\
		\midrule
		{SPINN--DEM} & $\underset{\pm29.70}{127.69}$ & $\underset{\pm209.73}{1291.20}$ & $\underset{\pm112.89}{760.34}$ & $\underset{\pm56.46}{646.89}$ & --- \\
		\bottomrule
	\end{tabular}\label{table7}
\caption{Thin-walled angle problem: Times to achieve ``engineering'' accuracy  ${\mathbb T}_{5\%}$ for the different values obtained by SPINN--DEM approach (in seconds). ``Engineering'' accuracy is defined by relative error which should be less or equal to 5$\%$ (${\mathbb L}_2 \leq 0.05$).}
\end{table}

%\begin{figure}[t!]
%	\centering
%	\includegraphics[width=1\textwidth]{figs/fig5_b.png}
%	\caption{(a) is prediction of a trained physics-informed neural network, (b) is reference solution, (c) are absolute values of difference between reference and predicted solutions, for the value $u_y$.}
%	\label{fig5b}
%\end{figure}
%
%\begin{figure}[t!]
%	\centering
%	\includegraphics[width=1\textwidth]{figs/fig5_c.png}
%	\caption{(a) is prediction of a trained physics-informed neural network, (b) is reference solution, (c) are absolute values of difference between reference and predicted solutions, for the value $u_z$.}
%	\label{fig5c}
%\end{figure}
%
%\begin{figure}[t!]
%	\centering
%	\includegraphics[width=1\textwidth]{figs/fig5_d.png}
%	\caption{(a) is prediction of a trained physics-informed neural network, (b) is reference solution, (c) are absolute values of difference between reference and predicted solutions, for the value $\sigma_{\rm vM}$.}
%	\label{fig5d}
%\end{figure}

	\newpage
	\section{Conclusion}
	
	In this article, we presented the method for solving elasticity problems based on separable physics-informed neural networks in conjunction with the deep energy method. The performance of this approach had been compared with the performances of the other method which are based on the solution of the partial differential equations which describe the elasticity problems. This comparison shows that the SPINN--DEM method has a convergence rate and accuracy that significantly exceeds the speed of convergence and accuracy of vanilla PINN and SPINN--PDE. {In terms of the speed of solving of the problems, this method approaches the classical numerical finite element methods (FEM)}. In addition, the authors' claim about original SPINN~\cite{cho2023} is confirmed, that SPINN-based solution methods consume much less GPU memory than the vanilla PINN. {Note another important feature of the methods of the PINN family is  a significantly smaller amount of memory occupied by a trained neural network compared to the stored data which generated by classical FEMs.} It was demonstrated that separable physics-informed neural networks in the framework of the deep energy method solve problems of the linear elasticity on complex geometries, which is not available to solve when using the PINNs in frames of partial differential equations. {The proposed method without significant changes can be used for a number of other elasticity problems~\cite{MENG2023116172,NING2023116430}.}

\newpage

\appendix

\section{Nomenclature}\label{appA}

The summarizes the main notations, abbreviations and symbols are given in table~8.
\begin{table}[h!]
	\centering
	\begin{tabular}{ll}
		\toprule
		Notation & Description  \\
		\midrule
		MLP & Multilayer perceptron  \\
		PDE & Partial differential equation\\
		PINN & Physic-informed neural network \\
		DEM & Deep energy method\\
		SPINN & Separable Physics-Informed Neural Networks\\
		SPINN--PDE & SPINN for the loss problem based on PDE\\
		SPINN--DEM & SPINN  in conjunction with the deep energy method (loss is energy loss)\\
		${\vec u}$ & solution of problem\\
		${\mathcal{B}}$ & boundary operator \\
		$\mathcal{R}$ & PDE residual \\
		${\vec u}^{\pmb \theta}$ & neural network representation of the problem solution \\
		${\pmb \theta}$ & vector of the trainable parameters of the neural network \\
		${\mathcal{L}(\pmb \theta)}$ & aggregate training loss\\
		\bottomrule
	\end{tabular}\label{tableAppA}
	\caption{Nomenclature}
\end{table}

\newpage
\bibliographystyle{unsrt} 
\bibliography{Eskin_Elasticity_Problem}

\end{document}